\documentclass[12pt,oneside,a4paper]{amsart}
\usepackage{geometry}
\usepackage[english]{babel}
\usepackage{amsfonts}
\usepackage{amsmath}
\usepackage{amssymb}
\usepackage{graphicx}
\usepackage{caption}
\usepackage{subcaption}
\usepackage{amsthm}
\usepackage{eurosym}
\usepackage[dvipsnames,usenames]{color}
\usepackage{cite}
\usepackage{soul}
\usepackage{tikz}
\usepackage{pgf}
\usepackage{lmodern} 
\usepackage[T1]{fontenc}
\title[Lower bounds]{\textbf{Some lower bounds for solutions of Schr\"odinger Evolutions}}
\author{Mikel Agirre}
\author{Luis Vega}
\subjclass[2010]{ 35Q41, 39A12}
\keywords{}
\date{\today}
\address{M. Aguirre Alonso: Basque Center for Applied Mathematics BCAM, Alameda de Mazarredo 14, 48009 Bilbao, Spain}
\email{mikelagi@outlook.com}
\address{L. Vega: Departamento de Matem\'aticas,  Universidad del Pa\'is Vasco UPV/EHU, apartado 644, 48080, Bilbao, Spain
\newline
\null\hspace{1.95cm}Basque Center for Applied Mathematics BCAM, Alameda de Mazarredo 14, 48009 Bilbao, Spain}
\email{luis.vega@ehu.eus}
\usepackage{hyperref} 
\newtheorem{theorem}{Theorem}[section]

\newtheorem{lemma}{Lemma}[section]

\theoremstyle{remark}

\newcommand{\supp}{\operatorname{supp}}

\begin{document}

\begin{abstract}
We present some lower bounds for regular solutions of Schr\"odinger equations with bounded and time dependent complex potentials. Assuming that the solution has  some positive mass at time zero within a ball of certain radius, we prove that this mass can be observed if one looks at the solution and its gradient in space-time regions outside of that ball.
\end{abstract}
\maketitle

\section{Introduction}

\noindent In this paper we are going to study the behavior of the solution of Schr\"odinger's initial value problem
\begin{equation}
\label{Schr1}\left\{ \begin{array}{ll}
\partial_t u = i(\Delta u + V(x,t)u) \\
u(x,0)=u_0(x)
\end{array} 
\right. 
\end{equation}

\noindent where $V(x,t)$ is a bounded  complex potential. For the case where the potential is identically zero, $V\equiv 0$, we can write down the explicit formula \begin{equation}
\label {Fourier}u(x,t)=e^{it\Delta}u_0(x)=\frac{1}{(2\pi)^{n/2}}\int_{\Bbb R^n} e^{-it|\xi|^2+ix\cdot\xi} \hat u_0(\xi)\,d\xi,\end{equation} 
\noindent where $e^{it\Delta}u_0(x)$ denotes the free solution and 
$$\hat f(\xi)= (2\pi)^{-n/2}\int_{\Bbb R^n} e^{-i\xi\cdot x}f(x) \,dx,$$

\noindent is the Fourier transform of $f$. Simple computations give that identity {\eqref{Fourier} can also be written as
\begin{align*}
\label {x,t}
u(x,t) & = \frac{1}{(4\pi it)^{n/2}}\int_{\Bbb R^n} e^{i \frac{|x-y|^2}{4t}}u_0(y)\,dy  \\
       & = \dfrac{e^{i\frac{|x|^2}{4t}}}{(4\pi it)^{n/2}}\int_{\Bbb R^n} e^{-i\frac{x}{2t}y}e^{i\frac{|y|^2}{4t}}u_0(y)\,dy.
\end{align*} 

\noindent As a consequence
\begin{equation*}
\label {Fourier2}e^{it\Delta}u_0(x)=\frac{e^{i|x|^2/4t}}{(2it)^{n/2}}\hat f_t\left(\frac{x}{2t}\right),
\end{equation*} 

\noindent with $f_t(x)=e^{i|x |^2/4t}u_0(x).$
This identity implies that to give size conditions of $u$ at two different times, say $t=0$ and $t=T$ is equivalent to give size conditions to $f_T$ and $\hat f_T$. This idea has been largely exploited by L. Escauriaza, C.E. Kenig, G. Ponce, and L. Vega, see for example \cite{BAMS} to revisit some classical Uncertainty Principles (UPs), as those by Hardy, Paley-Wiener, and Morgan. They give alternative proofs to these classical results using techniques of Partial Differential Equations, more concretely the so-called Carleman type inequalities. These UPs are rigidity results in the sense that the conclusions are that the only function satisfying the desired properties is either the trivial one, or some specific function, as for example the Gaussian in the case of Hardy's UP. Unlike this, the use of Carleman inequalities is rather flexible and allows perturbations by potentials $V\neq0$. Moreover, some lower bounds for general solutions are also proved.
 
\noindent All these results rely on assuming decay at two different times. The main purpose of this paper is to start to explore the way to reduce the hypothesis from two times to just one. Besides the fact that we consider this a very natural question by itself, our main motivation has been to try to adapt the PDE techniques to prove more sophisticated UPs as those proved by F. Nazarov in \cite{N}.


\noindent Our starting point is a general lower bound obtained in \cite{unique} (cf. Theorem 3.1) for solutions of \eqref{Schr1} for bounded potentials. The main assumption in that result is that the solution has some nontrivial mass  in a space-time cylinder with height $t\sim 1/2$ and the basis given by a ball of radius $1$. Then, the conclusion is that there exists a constant $c$ that depends just on the dimension such that the lower bound
\begin{equation}
\label{low}\int_0^1\int_{R-1<|x|<R}|u(x,t)|^2+|\nabla u(x,t)|^2dxdt\geq ce^{-cR^2},
\end{equation}

\noindent holds for all $R$ sufficiently large. In Theorem \ref{Schr} we will obtain lower bounds similar to \eqref{low} just assuming conditions at one time. In Theorem \ref{Uniq} we give some uniqueness results about solutions of \eqref{Schr1} for $x\in \Bbb R^n$.

\noindent The organization of this paper is as follows. In section 2 we gather the main results of the article. In sections 3 and 4 we give some auxilliary results necessary for the proof. In particular, we introduce the so called Appell's or pseudoconformal transformations. They give us the extra parameter we need to avoid to assume conditions at two different times. Here we use similar arguments to those in \cite{PW}. Finally the proof of Theorem \ref{Schr} is given in section 5.

\noindent {\bf{Acknowledgements:}} Both authors were partially supported  by an ERCEA Advanced Grant 2014 669689 - HADE,   by the MEIC projects MTM2014-53850-P and  SEV-2013-0323, and by the Basque Goverment BERC program.\\


\section{The main results}

\noindent By $\mathcal{H}^1_{loc}(\mathbb{R}^n)$ we mean the set of functions $f$ that together with their gradients are locally in $L^2(\Bbb R^n)$ . We have the following result.
\begin{theorem}\label{Schr}

Let $u\in\mathcal{C}([0,1]:\mathcal{H}^1_{loc}(\mathbb{R}^n))$ be a solution of
\[ \left\{ \begin{array}{ll}
\partial_t u = i(\Delta u + V(x,t)u) \\
u(x,0)=u_0(x),
\end{array} 
\right. \]

\noindent where $V\in L^{\infty}(\mathbb{R}^n\times[0,1])$ is a complex potential and
$$ \parallel V \parallel_{L^{\infty}(\mathbb{R}^n\times[0,1])}\leq L.$$
Let $R_0>0$ be such that for some $c_0>0$,
$$ \int_{B_{R_0}} \vert u_0 \vert^2 dx = c_0^2,$$
and let also $M\geq 4R_0+1$ so that
\begin{equation}
\sup_{0\leq t \leq 1} \int_{B_M} |u(x,t)|^2+|\nabla u(x,t)|^2dx = A^2 < +\infty.
\end{equation}

\noindent Then, there exist $t^*=\min\left( \frac{256A}{c_0L}, \  2^{-14}\left( \frac{c_0}{A}\right)^4, \ R_0^2, \  \frac{1}{L^2}\right)$ and a universal constant $c_n$ that depends just on the dimension such that if $0<t<t^*$,
\begin{equation}
\label{lw2}\frac{e^{c_n\frac{\rho^2}{t}}}{t}\int_{t/4}^{3t} \int_{||y|-\rho-\rho\frac{s}{t}|<4\rho\sqrt{t}}\vert u(y,s)\vert^2 + s\vert \nabla_y u(y,s) \vert^2\,dyds\geq c_0^2, \ \ \ \ \  R_0 \leq \rho \leq M.
\end{equation}

\end{theorem}




\noindent We state a few remarks about the theorem.

\noindent \textbf{Remark 1}: Observe that $M$ can be infinity.

\noindent \textbf{Remark 2}: 
\noindent  \textbf{(The periodic setting)} Take $\Bbb T=[0,2\pi].$ Assume now that we have a periodic in space solution $u\in\mathcal{C}([0,1];\mathcal{H}^1(\mathbb{T}))$ of the problem \eqref{Schr1} and  
$$ \int_{\mathbb{T}} \vert u_0 \vert^2 dx = c_0^2.$$ 

\noindent Choose $\rho = 2\pi$. Then the observability region has the form
$$\left| y \pm 2\pi\left(1+\frac{s}{t}\right)\right| < 8\pi\sqrt{t}, \ \ \ \ s\in[t/4,3t],$$
that by periodicity becomes
$$\left| y \pm 2\pi\frac{s}{t}\right | < 8\pi\sqrt{t}, \ \ \ \ s\in[t/4,3t].$$


\noindent Hence as a consequence of the above theorem there exists a universal constant $c_1$ such that
\begin{equation}
\label{lw3}\frac{e^{c_1/t}}{t}\int_{t/4}^{3t} \int_{|y\pm2\pi\frac{s}{t}|<8\pi\sqrt{t}}|u|^2+s|\nabla u|^2 dyds \geq c_0^2.
\end{equation}

\noindent This result can be seen as some kind of one dimensional observability inequality for periodic solutions of (1), as shown by N. Burq and M. Zworski in [3] for time independent potentials. Observe also that the variable $t$ is very small, as estated on Theorem \ref{Schr}. The drawback of \eqref {lw3} is that involves $\nabla u$, and it is a very natural question to know if this term is needed. Also, and because of this dependence on the gradient, it is not clear up to what extent \eqref {lw3} implies a controllability result. Nevertheless, observe that our result allows potentials that can depend on time and holds in any dimension. The geometry of  the controllability set becomes more complicated  as the dimension grows, and it would also be very natural to explore if, with the methods used in this paper, one could get closer to the results obtained for the two dimensional case by N. Anantharaman and F. Macia in [1], and J. Bourgain, N. Burq, and M.Zworski in [2].

\noindent Our second main result is a uniqueness one and it is an immediate consequence of Theorem \ref{Schr}, and  therefore the proof will be omitted. As far as we know this type of uniqueness result is completely new.
\begin{theorem}\label{Uniq}
Assume that for any $u_0\in H^1(\mathbb R^n)$ there exists a unique solution $u\in\mathcal C\big([0,1]\,:\,H^1(\mathbb R^n)\big)$ of
$$\left\{\begin{aligned}
\partial_t u&=i\left(\Delta+V(x,t)\right)u\quad x\in\mathbb R^n\,,\quad t\in (0,1)\\
u(x,0)&=u_0,
\end{aligned}
\right.$$
with $V\in L^{\infty}\big(\mathbb R^n\times [0,1]\big)$.
If $c_n$ is as in \eqref{lw2} and there exist $R_j$, $j\in\mathbb N$ such that  for all $j$
$$\lim_{t\downarrow 0}\displaystyle\frac 1t e^{c_n\frac{R_j^2}{t}}\int_{t/4}^{3t}\int_{||y| - R_j(1+s/t)| < 4R_j\sqrt{t}} 
|u(y,s)|^2+s\left|\nabla u(y,s)\right|^2 dy ds=0,$$
then $u\equiv 0$.
\end{theorem}

\noindent As a side result to Theorem \ref{Uniq} we could let the spatial parameter $\rho$ tend to infinity and obtain a similar conclusion.
\begin{theorem}\label{Uniq2}
Assuming the same conditions as in Theorem \ref{Uniq} and that there exists $t_0<t^*$ such that 
$$\lim_{\rho\rightarrow\infty}\displaystyle\frac 1{t_0} e^{c_n\frac{\rho^2}{t_0}}\int_{t_0/4}^{3t_0}\int_{||y| - \rho(1+s/t_0)| < 4\rho\sqrt{t_0}} 
|u(y,s)|^2+s\left|\nabla u(y,s)\right|^2 dy ds=0,$$
 then $u\equiv 0$.
\end{theorem}

\noindent \textbf{Remark 3}:
Notice that our results are perturbative and allow complex potentials that can depend on time. Therefore, it can be applied to solutions of non-linear equations as long as  a nice local in time well-posedness theory is available. We can proceed as done in \cite{BAMS} and consider for example
\begin{equation}
\label{NLS}
\partial_t u= i(\Delta u+ f(|u|^2)u),
\end{equation}
with $f: \Bbb R \rightarrow \Bbb R,$ and $f(0)=f'(0)=0$. Then, given two smooth solutions $u_1$ and $u_2$ of \eqref{NLS}, the difference $\omega=u_1-u_2$ satisfies an equation as \eqref{Schr1}, and therefore Theorems \ref{Uniq} and \ref{Uniq2} apply to $\omega$.


\section{Appell's conformal transformation and Carleman's estimate}

\noindent We use the following result from \cite{H} to generate a new family of solutions for Schr\"odinger's problem that depends on two parameters.

\begin{lemma}\label{Appell}
If $u(y,s)$ verifies
$$\partial_s u = i(\Delta u + V(y,s)u + F(y,s)), \ \ \ \ (y,s)\in\mathbb{R}^n\times [0,1]$$
and $\alpha$ and $\beta$ are positive, then 
$$ \tilde{u}(x,t)=\left(\frac{\sqrt{\alpha\beta}}{\alpha(1-t)+\beta t}\right)^{n/2}u\left(\frac{\sqrt{\alpha\beta}x}{\alpha(1-t)+\beta t}, \ \frac{\beta t}{\alpha(1-t)+\beta t}\right)e^{\frac{(\alpha - \beta)\vert x \vert^2}{4i(\alpha(1-t)+\beta t)}}$$
verifies
$$\partial_t \tilde{u} = i(\Delta \tilde{u} + \tilde{V}(x,t)\tilde{u} + \tilde{F}(x,t)), \ \ \ \ (x,t)\in\mathbb{R}^n\times [0,1]$$
with
$$\tilde{V}(x,t)=\frac{\alpha\beta}{(\alpha(1-t)+\beta t)^2}V\left(\frac{\sqrt{\alpha\beta}x}{\alpha(1-t)+\beta t}, \ \frac{\beta t}{\alpha(1-t)+\beta t}\right)$$
and
$$\tilde{F}(x,t)=\left(\frac{\sqrt{\alpha\beta}}{\alpha(1-t)+\beta t}\right)^{n/2+2}F\left(\frac{\sqrt{\alpha\beta}x}{\alpha(1-t)+\beta t}, \ \frac{\beta t}{\alpha(1-t)+\beta t}\right)e^{\frac{(\alpha - \beta)\vert x \vert^2}{4i(\alpha(1-t)+\beta t)}}.$$
\end{lemma}

\noindent Although the statement uses two parameters $\alpha$ and $\beta$ we are going to define $\gamma = \alpha / \beta$ and rewrite the dilations in a proper way. We need to be careful on how these functions alter the domains of integration in the proof. It is important to make a sensible use of the parameter $\gamma$ in relation to these functions. For this reason, we give some estimations in the next section. We will not be considering the function $F$ either, since we will not consider outer forces disturbing our system. 

\noindent Next we recall another result from \cite{unique}, that plays a fundamental role in the rest of the arguments. Therefore we include its proof for the sake of completeness.
\begin{lemma}[Carleman estimate]
Assume that $R>0$ and $\varphi:[0,1]\longrightarrow \mathbb{R}$ is a smooth function. Then, there exists $c_n=c(n, \ \parallel \varphi' \parallel_{\infty}+\parallel\varphi''\parallel_{\infty})>0$ such that the inequality 
\begin{equation}\label{carleman}
\frac{\sigma^{3/2}}{c_nR^2}\parallel e^{\sigma \left| \frac{x}{R}+\varphi(t)e_1\right|^2}g\parallel_2 \leq \parallel e^{\sigma \left| \frac{x}{R}+\varphi(t)e_1\right|^2}(i\partial_t+\Delta)g\parallel_2,
\end{equation}
holds when $\sigma \geq c_nR^2$ and $g\in\mathcal{C}_0^{\infty}(\mathbb{R}^{n+1})$ has its support contained in the set
$$\{(x,t) \ : \ |\frac{x}{R}+\varphi(t)e_1|\geq 1\}.$$
\end{lemma}

\noindent \textit{Proof.} Let $f=e^{\sigma |\frac{x}{R}+\varphi(t)e_1|^2}g$. Then,
$$e^{\sigma \left| \frac{x}{R}+\varphi(t)e_1\right|^2}(i\partial_t+\Delta)g = S_{\sigma}f-4\sigma A_{\alpha}f,$$

\noindent where 
$$S_{\sigma}=i\partial_t+\Delta + \frac{4\sigma^2}{R^2}|\frac{x}{R}+\varphi(t)e_1|^2,$$
$$A_{\sigma}= \frac{1}{R}(\frac{x}{R}+\varphi(t)e_1)\cdot\nabla + \frac{n}{2R^2}+\frac{i\varphi'}{2}(\frac{x_1}{R}+\varphi),$$

\noindent are the symmetric and anti-symmetric operators respectively. Thus,
$$A_{\sigma}^*=-A_{\sigma}, \ \ \ \ \ S_{\sigma}^*=S_{\sigma}$$

\noindent and
$$\parallel e^{\sigma \left| \frac{x}{R}+\varphi(t)e_1\right|^2}(i\partial_t+\Delta)g \parallel_2^2=\langle S_{\sigma}f-4\sigma A_{\sigma}f, \ S_{\sigma}f-4\sigma A_{\sigma}f\rangle$$
$$\geq -4\sigma \langle (S_{\sigma}A_{\sigma}-A_{\sigma}S_{\sigma})f,f\rangle = -4\sigma \langle [S_{\sigma},A_{\sigma}]f,f\rangle.$$

\noindent A calculation shows that 
$$[S_{\sigma},A_{\sigma}]= \frac{2}{R^2}\Delta - \frac{4\sigma^2}{R^4}\left|\frac{x}{R}+\varphi(t)e_1\right|^2-\frac{1}{2}\left[\left(\frac{x_1}{R}+\varphi(t)\right)\varphi''+\varphi'^2\right]+\frac{2i\varphi'}{R}\partial_{x_1},$$

\noindent and
$$\parallel e^{\sigma \left| \frac{x}{R}+\varphi(t)e_1\right|^2}(i\partial_t+\Delta)g \parallel_2^2 \geq \frac{16\sigma^3}{R^4}\int \left|\frac{x}{R}+\varphi e_1\right|^2\vert f \vert^2dxdt+\frac{8\sigma}{R^2}\int \vert \nabla f \vert^2dxdt$$
$$+2\sigma\int \left[\left(\frac{x_1}{R}+\varphi\right)\varphi''+\varphi'^2\right]\vert f \vert^2dxdt - \frac{8\sigma i}{R}\int \varphi'(\partial_{x_1}f) fdxdt.$$

\noindent Hence using the hypothesis on the support of $g$ and the Cauchy-Schwarz inequality, the absolute value of the last two terms can be bounded by a fraction of the first two terms on the right hand side when $\sigma \geq c_nR^2$ for some large $c_n$ depending on the dimension and $\parallel \varphi'\parallel_{\infty}+\parallel \varphi''\parallel_{\infty}$. Then the result follows.
\begin{flushright}
$\square$
\end{flushright}

\section{Some a priori estimates}

\noindent Before going into the proof of the theorem, we give some estimates of the support functions we are going to be using. At some point in the previous section we have talked about a transformation based on the conformal transformation and mentioned that the parameters $\alpha$ and $\beta$ will be replaced by $\gamma>0$ defined as the relation between the former ones, say $\alpha/ \beta$. We will also want this parameter $\gamma$ to be as big as possible so in principle we see it as $\gamma >> 1$. Having this in mind, we define the following functions:
\begin{equation}\label{20}
\alpha(t)= \frac{1}{\gamma^{1/2}(1-t)+\gamma^{-1/2}t}, 
\end{equation}
\begin{equation*}
s(t)=\frac{t}{\gamma(1-t)+t}, 
\end{equation*}
\begin{equation*}
\beta(t)=\frac{1}{1-t+\gamma^{-1}t}-\frac{1}{\gamma(1-t)+t}.
\end{equation*}

\noindent At some point on the proof there will be a change of variables so it is interesting to see how we can write $t$ in terms of $s$ and see how the measure changes with respect to $\gamma$. First we see that 
\begin{equation}\label{b}
t(s)=\frac{s\gamma}{1+s(\gamma-1)},
\end{equation}

\noindent and so
\begin{equation*}
dt=\frac{\gamma}{(1+s(\gamma-1))^2}ds.
\end{equation*}

\noindent Along the proof we encounter two different time intervals due to the definition of the cut-off functions. The first one is $[3/8,5/8]$. For this interval observe that $\alpha(t)$ can be estimated by
\begin{equation}\label{5}
\frac{1}{\gamma^{1/2}} \leq \alpha(t) \leq \frac{3}{\gamma^{1/2}},
\end{equation}

\noindent and the variable $s$ lives in
\begin{equation}\label{6}
I_s^1 = \left[\frac{3}{5\gamma+3},\frac{5}{3\gamma+5}\right]. 
\end{equation}

\noindent The length of this interval can be estimated by
\begin{equation}\label{7}
\frac{1}{4\gamma} \leq |I_s^1| \leq \frac{2}{\gamma},
\end{equation}

\noindent and the differential,
\begin{equation}\label{8}
\frac{\gamma}{8}ds \leq dt \leq \gamma ds.
\end{equation}

\noindent On the other hand, when $t\in[1/4,3/4]$ we can make the following estimations:
\begin{equation}\label{a}
\frac{1}{\gamma^{1/2}}\leq \alpha(t)\leq \frac{4}{\gamma^{1/2}} ,
\end{equation}
\begin{equation}
0\leq \beta(t)=\frac{1}{1-t+\gamma^{-1}t}-\frac{1}{\gamma(1-t)+t} \leq \frac{1}{1-t+\gamma^{-1}t} \leq 4.
\end{equation}

\noindent It should also be noticed that the variable $s$ lives on the interval
$$I_s^2 = \left[ \frac{1}{3\gamma+1}, \frac{3}{\gamma+3} \right],$$

\noindent and its length is bounded from above and below as follows
\begin{equation}\label{9}
\frac{1}{2\gamma} \leq |I_s^2| \leq \frac{3}{\gamma},
\end{equation}

\noindent which means that when $\gamma$ is large the variable $s$ has size $\gamma^{-1}$ and so we have the following estimation for the differential
\begin{equation}\label{10}
\frac{\gamma}{16}ds \leq dt \leq \gamma ds.
\end{equation}

\noindent Observe next that if we combine \eqref{20} with \eqref{b} we can write
\begin{equation}\label{d}
\sqrt{\gamma}\alpha(t(s)) = 1+s\gamma-s.
\end{equation}

\noindent Assume now that $u\in\mathcal{C}([0,1],\mathcal{H}^1_{loc}(\mathbb{R}^n))$ is a solution to \eqref{Schr1}. Then, the following identity holds
\begin{equation}\label{1}
|u(x,t)|^2-|u(x,0)|^2 = -2Im\int_0^t \left(div(u(x,s)\cdot \nabla \overline{u}(x,s)) + V(x,s) |u(x,s)|^2\right)ds,
\end{equation}

\noindent where $V$ is a complex bounded potential. The proof of this identity is the following: \\

\noindent First observe that 
$$div(u\nabla \bar{u})= |\nabla u|^2+u\Delta \bar{u} \Rightarrow u\Delta \bar{u}= div(u\nabla \bar{u})-|\nabla u |^2,$$
and so we compute the derivative on the second variable of the squared term 
\begin{align*}
\frac{d}{dt}|u(x,t)|^2 & = \frac{d}{dt}(u\bar{u}) = \partial_t u \bar{u}+ u \partial_t\bar{u} \\
							& = \overline{u\partial_t \bar{u}}+u\partial_t \bar{u}  = 2Re(u\partial_t\bar{u}) \\
							& = 2Re(u i(\Delta + V )\bar{u}) = 2Re(iu\Delta\bar{u}+iV|u|^2) \\
							& = -2Im(u\Delta\bar{u}+V|u|^2) \\
							& = -2Im(div(u\nabla\bar{u})-|\nabla u |^2+V|u|^2) \\
							& = -2Im(div(u\nabla\bar{u}) + V|u|^2),
\end{align*}

\noindent which concludes the proof.

\noindent We are ready now to discuss the proof of the main theorem on this paper.

\section{Proof of Theorem \ref{Schr}}

\noindent We follow very closely the arguments in \cite{PW}. The goal is to use the Carleman estimate \eqref{carleman} in a suitable way so that we can control both $u$ and $\nabla u$ by the initial data. For this purpose we want to build an auxilliary function $g$. First, let $\gamma$ be large enough, say $\gamma > 16$ and define $R=R_0\sqrt{\gamma}$. Define also the following cut-off functions, $\theta_R,\eta\in\mathcal{C}^{\infty}_0(\mathbb{R}^n)$, $\varphi\in\mathcal{C}^{\infty}([0,1])$

\[ \theta_R(x)= \left\{ \begin{array}{ll}
		1, & |x|\leq R \\
		0, & |x|\geq R+1
		\end{array}
	\right.
\ \ \ \ \  \eta(x)= \left\{ \begin{array}{ll}
		1, & |x|\geq 2 \\
		0, & |x|\leq 3/2
		\end{array}
	\right.
\]
\vspace{0.2cm}
\[
 \varphi(t)= \left\{ \begin{array}{ll}
		4, & t\in [3/8,5/8] \\
		0, & t\in [0,1/4]\cup [3/4,1].
		\end{array}
	\right.
\] 

\noindent For future purposes we will be assuming that $R\geq 2$. Next we use the conformal transformation (\ref{Appell}) on the solution $u$ to generate a new family of solutions depending on the parameter $\gamma$, say
$$v(x,t) = \alpha(t)^{n/2}u(\alpha(t)x,s(t))e^{-\frac{i}{4}\beta(t)\vert x\vert^2} \ \ \ , \ \ \ (x,t)\in\mathbb{R}^n\times [0,1]$$
where the functions $\alpha$, $\beta$ and $s$ were introduced on the previous section.

\noindent We use all the information gathered above to define the function $g$ as follows:
$$g(x,t) = \theta_R(x)\eta\left(\frac{x}{R}+\varphi(t)e_1\right)v(x,t) \ \ , \ \ (x,t)\in\mathbb{R}^n\times [0,1].$$
\noindent Observe that due to the nature of the test functions, $g$ is compactly supported and,

\noindent $\bullet$ $g=\theta_Rv$ on $(x,t)\in \{\vert x\vert\leq R+1\}\times [3/8,5/8]$ \\
$\bullet$ $\nabla_x v(x,t) = \alpha(t)^{n/2}e^{-\frac{i}{4}\beta(t)\vert x\vert^2}(\alpha(t)\nabla u - \frac{i}{2}\beta(t)x\cdot u)$ \\
$\bullet$ $\supp g \subseteq \{ \left| \frac{x}{R}+\varphi e_1 \right| \geq 1 \}$

\noindent where $u=u(\alpha(t)x,s(t))$. With the function we just defined, we are ready to use the Carleman estimate. Recall that for $\sigma \geq c_nR^2$ 
\begin{equation}
\label{C1} \dfrac{\sigma^{3/2}}{c_nR^2}\parallel e^{\sigma\left| \frac{x}{R}+\varphi(t)e_1\right|^2} g(x,t)\parallel_2 \  \leq \  \parallel e^{\sigma\left| \frac{x}{R}+\varphi(t)e_1\right|^2} (i\partial_t + \Delta)g(x,t)\parallel_2 .
\end{equation}

\noindent We need to work out both sides of the inequality. The goal is to give an estimation from below to the left hand side using the information we have about the initial data. Once this is done, we will find suitable upper estimates of the right hand side in order to hide the terms we don't need using the parameter $\sigma$. Let's thus look at the l.h.s. of the inequality above,
\begin{align*}
\parallel e^{\sigma\left| \frac{x}{R}+\varphi(t)e_1\right|^2} g(x,t)\parallel_2^2 & = \int_t\int_x e^{2\sigma\left| \frac{x}{R}+\varphi(t)e_1\right|^2}\vert g(x,t)\vert^2dxdt  \\
	& \geq e^{8\sigma}\int_{3/8}^{5/8}\int_{\vert x\vert\leq R+1} \vert \theta_R(x)v(x,t)\vert^2dxdt \\
	& = e^{8\sigma}\int_{3/8}^{5/8}\int_{\vert x\vert\leq R+1} \alpha(t)^n\vert \theta_R(x)u(\alpha(t)x,s(t))\vert^2dxdt \\
	& \geq e^{8\sigma}\frac{\gamma}{8}\int_{s(\frac{3}{5})}^{s(\frac{5}{3})}\int_{\vert y \vert \leq \alpha(t(s))(R+1)} \vert \theta_R(\alpha^{-1}(t(s))y)u(y,s)\vert^2dyds.
\end{align*}

\noindent We have made the change of variables $y=\alpha(t)x$ and $s=s(t)$ together with the estimate \eqref{8} on the differential and change of measure we mentioned on the previous section. Now we want to plug the initial data into the equation. To do so we measure the size of the difference between our function $u$ and the initial data $u_0$, say
$$ B= \left| \int_{s\sim \gamma^{-1}} \int_{|y| \leq \alpha(t(s))R} \theta_R^2(\alpha^{-1}(t(s))y)(|u|^2-|u_0|^2)dyds\right|.$$

\noindent Next we use \eqref{1} to obtain,
\begin{align*}
B & = \left| \int_{s\sim \gamma^{-1}} \int_y \theta_R^2 \left( -2Im \int_0^s(div(u\nabla \bar{u})+V|u|^2)ds' \right) dyds\right| \\
   & = \left| 2Im \int_{s\sim \gamma^{-1}}  \int_0^s\int_y \theta_R^2 (div(u\nabla \bar{u})+V|u|^2)dyds'ds\right| \\
   & \leq \left| 2Im \int_s\int_{s'}\int_y \theta_R^2 div(u\nabla \bar{u}) dyds'ds \right| + \left| 2Im \int_s\int_{s'}\int_y \theta_R^2 V|u|^2 dyds'ds \right| \\
   & = B_1 + B_2.
\end{align*}

\noindent Here we study the contribution of both integrals separately and see how to choose $\gamma$ in a suitable way depending on the parameters $c_0$, $A$ and $L$ so that we have a nice bound from below for the left hand side of Carleman's estimate on this particular case. For the estimation of both $B_1$ and $B_2$ we use \eqref{7} and \eqref{8} so that 
\begin{align*}
B_2 & = \left| 2Im \int_s\int_{s'}\int_y \theta_R^2(\alpha^{-1}(t(s))) V(y,s')|u(y,s')|^2 dyds'ds \right| \\
       & \leq \frac{4L}{\gamma} \int_{s'}\int_y |u|^2 dy ds' \\
       & \leq \frac{8L}{\gamma^2} \sup_{s'\sim\gamma^{-1}}\int_{|y|\leq 4R/\sqrt{\gamma}} |u|^2 dy \\
       & \leq \frac{8A^2L}{\gamma^2}.
\end{align*}

\noindent As for $B_1$ we have, using integration by parts,
\begin{align*}
B_1 & = \left| 4Im \int_{s\sim\gamma^{-1}}\int_0^s\alpha^{-1}(t(s)) \int_{|y|\leq \alpha(t(s))R} (\theta_R\nabla\theta_R)(u\nabla\bar{u})dyds'ds \right| \\
       & \leq \left| 4\sqrt{\gamma} Im \int_{s\sim\gamma^{-1}} \int_{s'\sim\gamma^{-1}} \int_y (\theta_R\nabla\theta_R)(u\nabla\bar{u})dyds'ds \right| \\
       & \leq 4\sqrt{\gamma} \int_{s\sim\gamma^{-1}} \int_{s'\sim\gamma^{-1}} \int_y |u\nabla\bar{u}| dy ds'ds \\
       & \leq 16\gamma^{-3/2} \sup_{s'\sim\gamma^{-1}}\int_{|y|\leq \alpha(t(s))R}  |u\nabla\bar{u}| dy  \\
       & \leq 8\gamma^{-3/2} \sup_{s'\sim\gamma^{-1}}\int_{|y|\leq 4R/\sqrt{\gamma}} |u|^2+|\nabla u|^2 dy \\
       & \leq \frac{8A^2}{\gamma^{3/2}}.
\end{align*}

\noindent Now if we put all together and remember that there was a factor $\gamma$ multiplying the equation, we have
$$\gamma B \leq \frac{8A^2}{\gamma^{1/2}}+ \frac{8A^2L}{\gamma} = \frac{8A^2}{\gamma^{1/2}}(1+L\gamma^{-1/2}).$$
Therefore, if we choose $\gamma \geq L^2$ we have that
$$\gamma B \leq \frac{16A^2}{\gamma^{1/2}}.$$

\noindent Now want to hide the contribution of $B$ using $c_0$. To do so we need to work out the extra term introduced when $B$ was defined, say
$$\gamma\int_s\int_y \theta_R^2|u_0|^2dyds.$$
\noindent Clearly, the inclusion $B_{\alpha(t(s))R}\subset B_{\alpha(t(s))(R+1)}$ and the definition of $\theta_R$ together with \eqref{7} gives us
\begin{align*}
\gamma\int_s\int_{B_{\alpha(t(s))(R+1)}} \theta_R^2 |u_0|^2dyds & \geq \gamma\int_s\int_{B_{\alpha(t(s))R}} \theta_R^2 |u_0|^2dyds \\
			& = \gamma\int_s\int_{B_{\alpha(t(s))R}}|u_0|^2dyds \\
			& \geq \gamma\frac{1}{4\gamma} \int_{B_{R_0}}|u_0|^2dy \\
			& \geq \frac{c_0^2}{4}.
\end{align*}

\noindent Thus if
$\frac{c_0^2}{8} \geq 16A^2\gamma^{-1/2}$, then
$$ \gamma \geq  2^{14}\left(\frac{A}{c_0} \right)^4,$$
\noindent and we can hide the contribution of $B$ inside $c_0^2$ and conclude
\begin{equation}\label{11}
\parallel e^{\sigma\left| \frac{x}{R}+\varphi(t)e_1\right|^2} g(x,t)\parallel_2^2 \geq \frac{e^{8\sigma}}{64}c_0^2.
\end{equation}

\noindent Now we study the right hand side of the Carleman estimate. First compute the operator to see how the supports of the resulting expressions change:
\begin{align*}
(i\partial_t+\Delta )g(x,t) & = \theta_R\eta \tilde{V} v +\theta_R(i\varphi'\partial_{x_1} \eta v + 2R^{-1}\nabla \eta \cdot\nabla v + R^{-2}\Delta \eta v) \\
& + \eta(2\nabla\theta_R \cdot \nabla v + \Delta\theta_R  v) \\
								   & = E_1+E_2+E_3,
\end{align*}
\noindent where $E_2$ and $E_3$ are supported in 

\noindent $\bullet$ $3/2\leq \left| \frac{x}{R}+\varphi(t)\right| \leq 2 \ \ , \ \ t\in[1/4,3/4] $ \\
\noindent $\bullet$ $\{R \leq \vert x \vert \leq R+1\} \times [1/4,3/4,]$ 

\noindent respectively. From the definition of the conformal transformation and the assumption on $V$ we have that $\vert\tilde{V}\vert \leq L\gamma^{-1} $ on the support of $g$. Thus,

\begin{align*}
\parallel e^{\sigma\left| \frac{x}{R}+\varphi(t)e_1\right|^2} (i\partial_t + \Delta)g(x,t)\parallel_2^2 &= \int_t\int_x e^{2\sigma\left| \frac{x}{R}+\varphi(t)e_1\right|^2}\vert(i\partial_t + \Delta) g(x,t)\vert^2dxdt  \\
								& \leq L^2\gamma^{-2}\parallel e^{\sigma\left| \frac{x}{R}+\varphi(t)e_1\right|^2}g\parallel_2^2 \\
								 & +e^{8\sigma}\int_{1/4}^{3/4}\int_{\vert x \vert \leq R+1} (\vert v \vert^2 + 4R^{-2}\vert \nabla v\vert^2)dxdt \\
								& + e^{72\sigma}\int_{1/4}^{3/4}\int_{R \leq \vert x \vert \leq R+1} (\vert v \vert^2 + \vert \nabla v\vert^2)dxdt \\
								& = L^2\gamma^{-2}\parallel e^{\sigma\left| \frac{x}{R}+\varphi(t)e_1\right|^2}g\parallel_2^2+ e^{8\sigma}I_1 +e^{72\sigma} I_2.					
\end{align*}

\noindent Observe that from \eqref{C1} the first term can be hidden on the left hand side of the inequality if 
\begin{equation}\label{2}
\frac{\sigma^{3/2}}{c_nR^2} \geq \frac{2L}{\gamma}
\end{equation}

\noindent So we only need to study the contribution of $I_1$ and $I_2$. To see things more clearly we split $I_1$ in the sub-integrals $I_{11}$, $I_{12}$, the first one measuring the contribution of $v$, and the second one doing the same for the gradient $\nabla v$.
\begin{align*}
I_{11} = \int_{1/4}^{3/4}\int_{\vert x\vert\leq R+1}\vert v\vert^2dxdt & = \int_{1/4}^{3/4}\int_{\vert x\vert\leq R+1}\alpha(t)^n\vert u(\alpha(t)x,s(t))\vert^2dxdt \\
				& \leq \gamma \int_{s(1/4)}^{s(3/4)}\int_{|y|\leq\alpha(t(s))(R+1)}\vert u(y,s)\vert^2dyds.
\end{align*}

\noindent Here we have simply made a change of variables and use the information we have about the behavior of the functions $\alpha(t)$ and $s(t)$ when $\gamma$ is big enough. As for $I_{12}$ we use the triangular inequality together with the change of variables $y=\alpha(t)x$ and the estimations with $\gamma$, as we see here
\begin{align*}
I_{12} & = \int_{1/4}^{3/4}\int_{\vert x\vert\leq R+1}4R^{-2}\vert \nabla v\vert^2dxdt \\
	   & = 4R^{-2}\int_{1/4}^{3/4}\int_{\vert x\vert\leq R+1}\alpha(t)^n\vert \alpha(t)\nabla u-\frac{i}{2}\beta(t)x u\vert^2dxdt \\
				& \leq 4R^{-2}\gamma\int_{s(1/4)}^{s(3/4)} \int_{|y|\leq \alpha(t(s))(R+1)}|\alpha(t(s))\nabla u(y,s) - \frac{i}{2}\beta(t(s))\alpha^{-1}(t(s))yu|^2dyds \\
				& \leq 4R^{-2}\gamma \int_{s\sim1/ \gamma} \int_{|y|\leq \alpha(t(s))(R+1)} \left( \frac{16}{\gamma}|\nabla u|^2 + 4(R+1)^2|u|^2 \right) dyds \\
				& =  4\gamma \int_{s\sim1/ \gamma} \int_{|y|\leq \alpha(t(s))(R+1)} \left( \frac{16}{R^2\gamma}|\nabla u|^2 + 4\left(1+\frac{1}{R}\right)^2|u|^2 \right) dyds.
\end{align*}

\noindent Using now that $R=R_0\sqrt{\gamma}\geq 2$ we get 
\begin{align*}
I_{12} & \leq 36\gamma \int_{s\sim\gamma^{-1}} \int_{|y|\leq \alpha(t(s))(R+1)} |u|^2 + \gamma^{-1}|\nabla u |^2dy ds. 
\end{align*}

\noindent And if we put both $I_{11}$ and $I_{12}$ together and use \eqref{9} we can estimate $I_1$ as follows
\begin{align*}
I_1 & \leq 72\gamma \int_{s\sim \frac{1}{\gamma}} \int_{|y|\leq \alpha(t(s))(R+1)} (|u|^2 + \gamma^{-1}|\nabla u|^2) dyds \\
      & \leq 216 \sup_{s\sim\gamma^{-1}} \int_{|y|\leq\alpha(t(s))(R+1)} (|u|^2 + \gamma^{-1}|\nabla u |^2) dy  \\
      & \leq 216 \sup_{s\sim\gamma^{-1}} \int_{|y|\leq M} (|u|^2 + \gamma^{-1}|\nabla u |^2) dy \\
      & \leq 216A^2 .
\end{align*}

\noindent  Following a similar computation we can estimate $I_2$ as,
$$I_2 \leq 32 \gamma^2 R_0^2 \int_{s\sim\frac{1}{\gamma}} \int_{\alpha(t(s))R\leq |y| \leq \alpha(t(s))(R+1)}|u|^2+\gamma^{-1}|\nabla u|^2 dyds.$$

\noindent Observe now that the spatial variable $y$ lives in a region of length $\alpha$. We would like to rewrite that region in terms of $\gamma$. Using (\ref{a}) and (\ref{d}) together with an appropriate estimation for $s$ we have that $I_2$ can be written as 
$$I_2 \leq \gamma^2 32R_0^2\int_{s\sim \gamma^{-1}}  \int_{||y|-R_0-R_0s\gamma|<\frac{4R_0}{\sqrt{\gamma}}} \vert u\vert^2 + \gamma^{-1}\vert \nabla u\vert^2dyds. $$

\noindent  Now if we put everything together we have the following inequality
\begin{align}\label{4}
\frac{\sigma^{3/2}}{c_nR^2}\frac{c_0}{16}  \leq 16 A + 6\gamma R_0 e^{36\sigma} \left(\int_{s\sim\gamma^{-1}} \int_{||y|-R_0-R_0s\gamma|<\frac{4R_0}{\sqrt{\gamma}}}( \vert u\vert^2 + \gamma^{-1}\vert \nabla u\vert^2)dyds \right)^{1/2}.
\end{align}

\noindent To hide the first term of the right hand side inside the left hand side we ask the following 
\begin{equation}\label{3}
\frac{\sigma^{3/2}}{c_nR^2}  \geq \frac{512 A}{c_0} \Longrightarrow  \sigma \geq \hat{c_n}R^{4/3},
\end{equation}
for some universal $\hat{c_n}$ that depends on $c_n$, $c_0$ and $A$.
\noindent Since we want \eqref{2} and \eqref{3} to be satisfied we impose the following condition on the parameter $\gamma$
\begin{equation}
\frac{2L}{\gamma}\leq \frac{512 A}{c_0} \Longrightarrow \gamma \geq \frac{c_0L}{256A}.
\end{equation}

\noindent And thus, whenever 
$$\gamma \geq \max \left( \frac{c_0L}{256A}, \ 2^{14}\left(\frac{A}{c_0}\right)^4, \ \frac{1}{R_0^2}, \ L^2 \right), $$

\noindent we can hide the contribution of the first term on the right hand side of \eqref{4} into the left hand side, so
$$ \frac{\sigma^{3/2}}{c_nR^2} \frac{c_0}{32} \leq 16R_0\gamma e^{36\sigma}  \left(\int_{s\sim\gamma^{-1}} \int_{||y|-R_0-R_0s\gamma|<\frac{4R_0}{\sqrt{\gamma}}}( \vert u\vert^2 + \gamma^{-1}\vert \nabla u\vert^2)dyds \right)^{1/2}.$$

\noindent On the other hand $\sigma$ has to be greater than $c_nR^2$ according to Carleman's estimate, which is a stronger condition than the one we just found. Hence if $\sigma = 64c_nR^2$ we have that for some universal constant $c_{n,1}$ which depends only on $c_n$,
$$c_0 \leq \left(\gamma e^{c_{n,1}R_0^2\gamma}\int_{s\sim \gamma^{-1}}\int_{||y|-R_0-R_0s\gamma|<\frac{4R_0}{\sqrt{\gamma}}} (|u|^2+s|\nabla u |^2)dyds \right)^{1/2}, $$

\noindent now if we rename $\gamma^{-1} \equiv t$ observe that 
$$t \leq \min \left( \frac{256A}{c_0L}, \  2^{-14}\left( \frac{c_0}{A}\right)^4, \ R_0^2, \  \frac{1}{L^2} \right) = t^*,$$
and
$$c_0^2\leq \frac{e^{c_{n,1}\frac{R_0^2}{t}}}{t} \int_{t/4}^{3t} \int_{||y|-R_0-R_0\frac{s}{t}|<4R_0\sqrt{t}} (|u|^2+s|\nabla u|^2)dyds,$$

\noindent as we wanted to see.

\begin{flushright}
$\square$
\end{flushright}

\noindent On the statement of the theorem we write $c_n$ for simplicity. Observe also that if we take $\rho\in [R_0,M]$ the result will still be true. This happens because no matter what $\rho$ we choose, $c_0$ does not change.

\end{document}